\documentclass[12pt]{amsart}
\usepackage{ amsmath, amsthm, amsfonts, amssymb, color, xcolor}
 \usepackage{mathrsfs}
\usepackage{amsfonts, amsmath}
 \usepackage{amsmath,amstext,amsthm,amssymb,amsxtra,verbatim}

\begin{document}

\newfam\msbfam
\font\tenmsb=msbm10    \textfont\msbfam=\tenmsb \font\sevenmsb=msbm7
\scriptfont\msbfam=\sevenmsb \font\fivemsb=msbm5
\scriptscriptfont\msbfam=\fivemsb
\def\Bbb{\fam\msbfam \tenmsb}
\newfam\bigfam
\font\tenbig=msbm10 scaled \magstep2   \textfont\bigfam=\tenbig
\font\sevenbig=msbm7 scaled \magstep2 \scriptfont\bigfam=\sevenbig
\font\fivebig=msbm5 scaled \magstep2

\scriptscriptfont\bigfam=\fivebig

\def\dsum{\displaystyle\sum}
\def\dfrac{\displaystyle\frac}
\def\dlim{\displaystyle\lim}
\def\dint{\displaystyle\int}
\def\dsup{\displaystyle\sup}
\def\nln{\newline}

\newtheorem{thm}{Theorem}
\newtheorem{lem}{Lemma}[section]
\newtheorem{prop}{Proposition}[section]
\newtheorem{rem}{Remark}
\newtheorem{cor}{Corollary}[section]
\newtheorem{defn}{Definition}[section]
\newtheorem{pf}{\it Proof.}
\renewcommand\thepf{}

\newcommand{\beq}{\begin{equation}}
\newcommand{\eeq}{\end{equation}}

\title[ On the  norm of the operator $aI+bH$]
{On the norm of the operator $aI+bH$ on $L^p(\mathbb R)$}

\author{Yong Ding}

\address{Yong Ding, Laboratory of Mathematics and Complex Systems, School of Mathematical Sciences, Beijing Normal University, Ministry of Education of China, Beijing 100875, China }
\email{dingy@bnu.edu.cn }

\author{Loukas Grafakos}

\address{Loukas Grafakos, Department of Mathematics, University of Missouri, Columbia MO 65211, USA}
\email{grafakosl@missouri.edu}

\author{Kai Zhu${}^1$}

\address{Kai Zhu, School of Mathematical Sciences, Beijing Normal University, Beijing 100875, China}
\email{kaizhu0116@126.com }

\thanks{{\it Keywords and phases:} Best constants, Hilbert transform.}
\thanks{The first author is supported by NSFC(11371057,11471033,11571160),SRFDP(20130003110003) and the Fundamental Research Funds for the Central Universities(2014KJJCA10).}
\thanks{The second author  would like to acknowledge the support of Simons Foundation.}
\thanks{The third author would like to thank the China Scholarship Council for its support.}
\thanks{${}^1$Corresponding author.}

\begin{abstract}
We provide a direct proof of the  following theorem of Kalton, Hollenbeck, and Verbitsky \cite{HKV}: let $H$ be
the Hilbert transform  and let $a,b$ be real constants. Then for $1<p<\infty$ the
norm of the operator   $aI+bH$ from $L^p(\mathbb R)$ to $L^p(\mathbb R)$  is equal to
$$
\bigg(\max_{x\in \Bbb R}\frac{|ax-b+(bx+a)\tan \frac{\pi}{2p}|^p+|ax-b-(bx+a)\tan \frac{\pi}{2p}|^p}{|x+\tan \frac{\pi}{2p}|^p+|x-\tan \frac{\pi}{2p}|^p} \bigg)^{\frac 1p}.
$$
Our proof avoids passing through the analogous result for the conjugate function on the circle, as in \cite{HKV}, and is given directly on the line.   We also provide new approximate
extremals for $aI+bH$ in the case $p>2$.
\end{abstract}

\maketitle

\section{Introduction}

In this note we revisit the celebrated result of Kalton, Hollenbeck, and Verbitsky \cite{HKV} concerning the value of the
norm of the operator $aI+bH$ from $L^p(\mathbb R)$ to $L^p(\mathbb R)$ for $1<p<\infty$ and $a, b$   real constants. We provide a self-contained direct proof of this result on the real line. The original proof in \cite{HKV} was given for the conjugate function on the circle in lieu of the Hilbert transform and the corresponding  result for the line  was obtained from the periodic case via a transference-type argument due to Zygmund~\cite[Ch XVI, Th. 3.8]{Z} known as ``blowing up the circle". Here we work directly with the Hilbert transform on the line, using an idea contained in \cite{G2} and \cite{GS}, which is based on applying subharmonicity on the boundary of a suitable family of discs that fill up the upper half space  as their radii tend to infinity. The main estimates needed for our proof (Lemmas~\ref{Lemma 2.1} and \ref{Lemma 2.2}) are as in   \cite{HKV} but are included in this note for the sake of completeness (with a minor adjustment). The new contributions of this article are contained in Sections~\ref{Sec4} and \ref{Sec5}.  In Section~\ref{Sec4} we use a limiting argument and  subharmonicity  to prove the claimed bound for $aI+bH$. We obtain the approximate extremals for the operators $aI+bH$ in Section~\ref{Sec5}; these are based on these for the Hilbert transform which first appeared in Gohberg and Krupnik~\cite{GK} for $1<p<2$ and were also used by Pichorides~\cite{Pi}. We find new approximate extremals for the Hilbert transform for $2<p<\infty$ in Section~\ref{Sec5} and we use them to construct corresponding approximate extremals for $aI+bH$ for this range of $p$'s.
We note that the case $a=0, b=1$ of this result was proved by Pichorides \cite{Pi} and B. Cole (unpublished, see
\cite{Gamelin}), while the case  $a=0$, $b=1$, $p=2^m$, $m=1,2,\dots$, was  obtained four years  earlier by
Gohberg and Krupnik~\cite{GK}. For a short history on this topic we refer to Laeng~\cite{Laeng}.
It is noteworthy that the operator norm  of the Hilbert transform on $L^p$ is also the norm of other operators, for
instance of the segment multiplier; on this see De Carli and Laeng
\cite{decarli-laeng}.

\section{The norm of $aI+bH$}
Denote the identity operator by $I$. The Hilbert transform on the real line is defined by
$$
Hf(x)=p.v.\frac{1}{\pi}\int_{\Bbb R}\frac{f(t)}{x-t}\, dt
$$
for a smooth function with compact support.
For $a,b \in \Bbb R$, define
\beq\label{(3)}
B_p=\max_{x\in \Bbb R}\frac{|ax-b+(bx+a)\tan \gamma|^p+|ax-b-(bx+a)\tan \gamma|^p}{|x+\tan \gamma|^p+|x-\tan \gamma|^p},
\eeq
where $\gamma =\frac{\pi}{2p}$. $B_p$ can be defined equivalently by
\beq\label{(4)}
B_p=(a^2+b^2)^{p/2}\max_{0\leq \theta \leq 2\pi}\frac{|\cos(\theta +\theta_0)|^p+|\cos(\theta +\theta_0+\frac{\pi}{p})|^p}{|\cos \theta|^p+|\cos(\theta +\frac{\pi}{p})|^p},
\eeq
where $\tan \theta_0=b/a$. By letting $\theta=-\vartheta-\pi/p$,
\beq\label{(New1.0)}
B_p=(a^2+b^2)^{p/2}\max_{0\leq \vartheta \leq 2\pi}\frac{|\cos(\vartheta -\theta_0)|^p+|\cos(\vartheta -\theta_0+\frac{\pi}{p})|^p}{|\cos \vartheta|^p+|\cos(\vartheta +\frac{\pi}{p})|^p}.
\eeq
Our goal is to provide a proof of the following result in \cite{HKV}:

\begin{thm}\label{T1.1}\cite{HKV} Let $1<p<\infty$ and for $a,b \in \Bbb R$. Then
for all smooth functions with compact support $f$ on the line we have
$$
\|(aI+bH)f\|^p_{L^p({\Bbb R})}\leq B_p\|f\|^p_{L^p({\Bbb R})}
$$
where the constant $B_p$ is sharp. In other words,
$$
\|aI+bH\|_{L^p({\Bbb R})\to L^p({\Bbb R}) }= B_p^{\frac{1}{p}}  .
$$
\end{thm}

Without lose of generality, we assume that $a=\cos \theta_0, b=\sin \theta_0$, so that $a^2+b^2=1$.
As $aI+bH$ maps real-valued functions to real-valued functions,
in view of the Marcinkiewicz and Zygmund theorem~\cite{MZ} (see also~\cite[Theorem 5.5.1]{G}), the norm of $aI+bH$ on $L^p(\Bbb R) $  and on $ L^p({\Bbb C}) $    are equal.\footnote{for operators that do not map real-valued functions to real-valued functions, these norms may not be equal; for instance this is the case for the Riesz projections, see \cite{HV}.} Thus we may work with a nice real-valued function $f$ in the proof of Theorem~\ref{T1.1}.

\section{Some Lemmas}
In this section we provide two auxiliary results that are crucial in the proof of the main theorem.

\begin{lem}\label{Lemma 2.1}\cite{HKV} Suppose $p>1/2$, $p\neq 1$, and $F$ is a p-homogeneous continuous function on $\Bbb C$. Suppose there is a sector $S$ so that $F$ is subharmonic on $S$ and superharmonic on the complementary sector $S'$. Suppose further there is no nontrivial sector on which $F$ is harmonic. Suppose that $F(z)+F(e^{i\pi /p}z)\geq 0$ for all $z$, and there exists $z_0\neq 0$ so that $F(z_0)+F(e^{i\pi /p}z_0)=0$. Then there is a continuous p-homogeneous subharmonic function $G$ with $G(z)\leq F(z)$ for all $z\in \Bbb C$.
\end{lem}

\emph{Proof.} Lemma~\ref{Lemma 2.1} is  a restatement of  Theorem 3.5 in \cite{HKV}. We only provide a sketch below     making a minor  modification in the proof in  \cite{HKV} (i.e.,   definition of $h$ in \eqref{defh}).

We can suppose there exists $z_0$ with $|z_0|=1$ so that $F(z_0)+F(e^{i\pi /p}z_0)=0$. Let $z_0=e^{it_0}$, $z_1=e^{i(t_0+\pi /p)}$, since $p>1/2$, there exists $\epsilon >0$ such that $t_0-\epsilon <t_0<t_0+\pi/p<t_0+2\pi-\epsilon $. Write $F(re^{it})=r^pf(pt)$, where $f$ is a $2p\pi$-periodic function on $\Bbb R$. By Proposition 3.3 in \cite{HKV}, if $I$ is any interval so that $e^{ix/p}\in S$ for $x\in I$, then $f$ is trigonometrically convex on $I$, and if $e^{ix/p}\in S'$ for $x\in I$, then $f$ is trigonometrically concave on $I$. At least one of $z_0,z_1$ is contained in  $ S $; let us suppose that $z_0\in S$.  The function $f(x)+f(x+\pi)$ has minimum at $pt_0$, hence
$
f'_-(pt_0)+f'_-(pt_0+\pi)\leq 0,
$
$
f'_+(pt_0)+f'_+(pt_0+\pi)\geq 0.
$
This implies that there exist $a$ and $b$ such that $a+b=0$ and $f'_-(pt_0)\leq a\leq f'_+(pt_0)$ and $f'_-(pt_0+\pi)\leq b\leq f'_+(pt_0+\pi)$. Now define
\beq\label{defh}
h(x)=f(pt_0)\cos (x-pt_0)+a\sin (x-pt_0).
\eeq
Then by Lemma 3.1 in \cite{HKV}, $h \leq f $ on a neighborhood of $pt_0$.   Lemma 3.2 in \cite{HKV} implies that $h(x)\leq f(x)$ for $p\alpha +2p\pi \leq x <pt_0+\pi +\delta$ and for $pt_0\leq x\leq p\beta$. By the Phragm\'{e}n-Lindel\"{o}f theorem  (\cite{L}) we
obtain that $h \leq f $   in a neighborhood of $[pt_0,pt_0+\pi]$.

Let $T=\{re^{i\theta}:r>0,t_0<\theta <t_0+\frac{\pi}{p} \},$ define $H(re^{i\theta})=r^ph(p\theta)$ for $t_0<\theta<t_0+\frac{\pi}{p}$ and
$$
G(z)=\begin{cases}
H(z)\qquad \text{if  $z\in T$},\\
F(z)\qquad \text{if  $z\notin T$}.
\end{cases}
$$
Then $G(z)\leq F(z)$ for all $z\in \Bbb C$, by this we mean $\{re^{it}: r>0, t_0-\epsilon \leq t<t_0+2\pi-\epsilon \}$, and $G$ is subharmonic on both T and its complementary sector $T'$. It is easy to see $G$ is then subharmonic on ${\Bbb C}\backslash \{0\}$ since $h \leq f $ in a neighborhood of $pt_0$ and $pt_0+\pi$. Finally $h(x)+h(x+\pi)=0$ and Lemma 3.1 in \cite{HKV} imply that $G(z)+G(e^{i\pi /p}z)\geq 0$ for all $z$. Integrating   over a circle around $0$ yields the subharmonicity of $G$ at $0$.
\qed

\medskip
Next we have a version of Lemma  4.2 in \cite{HKV}  in which we provide an explicit formula for the
subharmonic  function $G$.

\begin{lem}\label{Lemma 2.2} Let $1<p<\infty$, $B_p$ be given by \eqref{(3)}, $T=\{re^{it}:r>0,t_0<t<t_0+\frac{\pi}{p}\}$, where $t_0$ is the value that makes right part of \eqref{(New1.0)} attain its maximum, and there exists $\varepsilon>0$ such that $t_0-\varepsilon<t_0<t_0+\pi/p<t_0+\pi-\varepsilon$. Let $z=re^{it}, z_0=re^{it_0}, G(z)=G(re^{it})$ be $\pi$-periodic of $t$ and when $t_0-\varepsilon<t<t_0+\pi-\varepsilon$:
$$
G(z)=\begin{cases}
B_p|\mathrm{Re}z_0|^{p-1}\textup{sgn}(\mathrm{Re}z_0)\mathrm{Re}[(\frac{z}{z_0})^pz_0]-|a\mathrm{Re}z_0+b\mathrm{Im}z_0|^{p-1}\\
\quad \quad\times \textup{sgn}(a\mathrm{Re}z_0+b\mathrm{Im}z_0)(a\mathrm{Re}[(\frac{z}{z_0})^pz_0]+b\mathrm{Im}[(\frac{z}{z_0})^pz_0]), &\textup{if} z\in T  \\
B_p|\mathrm{Re}z|^p-|a\mathrm{Re}z+b\mathrm{Im}z|^p, \quad \quad \quad\quad\quad\quad\quad\quad \,\,\qquad &\textup{if} z \notin T.
\end{cases}
$$
Then $G(z)$ is subharmonic on $\Bbb C$ and satisfies
\beq\label{(5)}
|a\mathrm{Re}z+b\mathrm{Im}z|^p\leq B_p|\mathrm{Re}z|^p-G(z).
\eeq
for all $z\in \Bbb C$.
\end{lem}

\emph{Proof.} The case $b=0$ is trivial, so we assume $b\neq 0$, and we may further assume that $a^2+b^2=1$. Let $F(z)=B_p|\mathrm{Re}z|^p-|a\mathrm{Re}z+b\mathrm{Im}z|^p$. Then  $F(re^{it})=r^pf(t)$, where $f(t)=B_p|\cos t|^p-|a\cos t+b\sin t|^p$ is $\pi$-periodic and continuously differentiable. The definition in \eqref{(4)} implies that
$$
\min_{0\leq t\leq 2\pi}[f(t)+f(t+\pi /p)]=0.
$$
We observe that $\Delta F\geq 0$ is equivalent to
$$
B_p|\mathrm{Re}z|^{p-2}\geq |a\mathrm{Re}z+b\mathrm{Im}z|^{p-2}.
$$
In order for $F(z)$ to be subharmonic, the following must be true:
$$
|a+b\tan t|^{p-2}\leq B_p.
$$
We can see that for $p\neq 2$ there will be two separate ``double sectors'' where $F(z)$ is subharmonic, and superharmonic in their complement. So let $\widetilde{p}=p/2, \widetilde{t_0}=2t_0$, define $\widetilde{F}(z)=F(z^{1/2})$, then $\widetilde{F}$ is $\widetilde{p}$-homogeneous and satisfies the hypotheses of Lemma~\ref{Lemma 2.1} with $\widetilde{p}$ and $\widetilde{t_0}$. Write $\widetilde{F}(re^{it})=r^{\widetilde{p}}\widetilde{f}(\widetilde{p}t)$, where
$$
\widetilde{f}(t)=B_p|\cos(t/p)|^p-|a\cos(t/p)+b\sin(t/p)|^p.
$$
 We can get
\beq\label{(6)}
\widetilde{f}(\widetilde{p}\widetilde{t_0})=B_p|\cos t_0|^p-|\cos(t_0-\theta_0)|^p,
\eeq
\beq\label{(7)}
\widetilde{f}'_-(\widetilde{p}\widetilde{t_0})=\widetilde{f}'_+(\widetilde{p}\widetilde{t_0})=-B_p\frac{|\cos t_0|^p}{\cos t_0}\sin t_0+\frac{|\cos (t_0-\theta_0)|^p}{\cos (t_0-\theta_0)}\sin (t_0-\theta_0),
\eeq
where $\tan \theta_0=b/a$. By the proof of Lemma~\ref{Lemma 2.1}, let
$$
h(x)=\widetilde{f}(\widetilde{p}\widetilde{t_0})\cos (x-\widetilde{p}\widetilde{t_0})+\widetilde{f}'_+(\widetilde{p}\widetilde{t_0})\sin (x-\widetilde{p}\widetilde{t_0}),
$$
then $h(x)\leq \widetilde{f}(x)$ for all $x$ in a neighborhood of $[\widetilde{p}\widetilde{t_0},\widetilde{p}\widetilde{t_0}+\pi]$.

Let $\widetilde{T}=\{re^{it}:r>0,\widetilde{t_0}<t <\widetilde{t_0}+\frac{\pi}{\widetilde{p}} \}$, and $H(re^{it})=r^{\widetilde{p}}h(\widetilde{p}t)$ for $\widetilde{t_0}<t<\widetilde{t_0}+\frac{\pi}{\widetilde{p}}$, let
$$
\widetilde{G}(z)=\begin{cases}
H(z)=H(re^{it})=r^{\widetilde{p}}[\widetilde{f}(\widetilde{p}\widetilde{t_0})\cos(\widetilde{p}t-\widetilde{p}\widetilde{t_0})+\widetilde{f}'_
+(\widetilde{p}\widetilde{t_0})\sin(\widetilde{p}t-\widetilde{p}\widetilde{t_0})]  \quad \text{if } z\in \widetilde{T},\\
\widetilde{F}(z)=r^{\widetilde{p}}(B_p|\cos \frac{t}{2}|^p-|a\cos \frac{t}{2}+b\sin \frac{t}{2}|^p)  \qquad\qquad\qquad\qquad\qquad\:\: \text{if } z\notin \widetilde{T}.
\end{cases}
$$
So let $\epsilon=2\varepsilon$, $\widetilde{G}$ is subharmonic and $\widetilde{G}(z)\leq \widetilde{F}(z)$ on $\{re^{it}: r>0, \widetilde{t_0}-\epsilon \leq t<\widetilde{t_0}+2\pi-\epsilon \}$ by Lemma~\ref{Lemma 2.1}. Now let $G(z)=\widetilde{G}(z^2)$, clearly $G$ is $p$-homogeneous and satisfies $G(z)\leq F(z)$ for $\{re^{it}: r>0, t_0-\varepsilon<t<t_0+\pi-\varepsilon \}$.
Since $z^2$ is holomorphic, $G(z)$ is also subharmonic on $\{re^{it}: r>0, t_0-\varepsilon<t<t_0+\pi-\varepsilon \}$.
Now let function $G(z)=G(re^{it})$ be $\pi$-periodic. For $t_0-\varepsilon\leq t<t_0+\pi-\varepsilon$, by \eqref{(6)}, \eqref{(7)} and $G(z)=\widetilde{G}(z^2)$ we have:
$$
G(z)=\begin{cases}
r^p[B_p \dfrac{|\cos t_0|^p}{\cos t_0}\cos(p(t-t_0)+t_0)-\dfrac{|\cos (t_0-\theta_0)|^p}{\cos (t_0-\theta_0)}&\cos(p(t-t_0)+t_0-\theta_0)],\\
&\text{if}\, z\in T ,\\
r^p(B_p|\cos t|^p-|\cos(t-\theta_0)|^p),  &\text{if}\, z \notin T,
\end{cases}
$$
where $\tan \theta_0=b/a$. It is easy to see $G(z)\leq F(z)$ for all $z\in \Bbb C$, by this we mean $\{re^{it}: r>0, t_0-\varepsilon \leq t<t_0+2\pi-\varepsilon \}$, so we get \eqref{(5)}.
Using similar proof as Lemma~\ref{Lemma 2.1} and the periodicity of $G$, we can get $G(z)$ is also subharmonic on $\Bbb C$.
Since $z_0=re^{it_0}$, the above formula is equivalent to
$$
G(z)=\begin{cases}
B_p|\mathrm{Re}z_0|^{p-1}\textup{sgn}(\mathrm{Re}z_0)\mathrm{Re}[(\frac{z}{z_0})^pz_0]-|a\mathrm{Re}z_0+b\mathrm{Im}z_0|^{p-1}\\
\quad \times \textup{sgn}(a\mathrm{Re}z_0+b\mathrm{Im}z_0)(a\mathrm{Re}[(\frac{z}{z_0})^pz_0]+b\mathrm{Im}[(\frac{z}{z_0})^pz_0]), & \hspace{-.1in} \text{if}\, z\in T , \\
B_p|\mathrm{Re}z|^p-|a\mathrm{Re}z+b\mathrm{Im}z|^p,  & \hspace{-.1in}\text{if}\, z \notin T.
\end{cases}
$$
This completes the proof of the lemma. \qed

\section{Proof of Theorem~\ref{T1.1}}\label{Sec4}
If $p=2$, then obviously
$$
\|aI+bH\|^2_{L^2({\Bbb R})\rightarrow L^2({\Bbb R})}=a^2+b^2=B_2,
$$
so we can assume $p\neq 2$. Consider the holomorphic extension of $f(x)+iH(f)(x)$ on the upper half space given by
$$
u(z)+iv(z)=\frac{i}{\pi}\int^{+\infty}_{-\infty}\frac{f(t)}{z-t}dt,\quad u,v \;\text{real-valued.}
$$
Let $G(z)$ be given by Lemma~\ref{Lemma 2.2}, our next step is to use Lemma~\ref{Lemma 2.2} and replace $z$ with $h(z)=u(z)+iv(z)$. Since $h(z)$ is holomorphic and $G$ is subharmonic, it follows that $G(h(z))$ is subharmonic on the upper half space. We note that (\cite{G2})
$$
|u(x+iy)|+|v(x+iy)|\leq \frac{C_f}{1+|x|+|y|}.
$$
By Lemma~\ref{Lemma 2.2}, we have that $|G(z) | \le C |z|^p$, hence
$$
|G(h(z))|\leq C|h(z )|^p \leq C(|u(z)|+|v(z)|)^p.
$$
So
\beq\label{(9)}
|G(h(z))|\leq \frac{C^p_f}{(1+|x|+|y|)^p}
\eeq
where $z=x+iy$. The boundary values of $G(h(z))$ are $G(h(x+i0))$.

The following part of the argument is based on \cite{GS}.
For $R>100$, consider the circle with center $(0,R)$ and radius $R'=R-R^{-1}$, denote by
$$
C^U_R=\{iR+R'e^{i\phi}:-\pi/4\leq \phi \leq 5\pi/4\}
$$
and
$$
C^L_R=\{iR+R'e^{i\phi}:5\pi/4\leq \phi \leq 7\pi/4\}.
$$
It follows from the subharmonicity of $G(h(z))$ that
\beq\label{(10)}
\int_{C^U_R}G(h(z))ds+\int_{C^L_R}G(h(z))ds\geq 2\pi R'G(h(iR)).
\eeq
Clearly \eqref{(9)} implies that
\beq\label{(11)}
|R'G(h(iR))|\leq R'\dfrac{C}{(1+R)^p}\rightarrow 0 \qquad \text{as}\; R\rightarrow \infty ,
\eeq
and that
\beq\label{(12)}
\bigg|\int_{C^U_R}G(h(z))ds  \bigg|\leq R'\dfrac{C}{(1+R)^p}\rightarrow 0 \qquad \text{as}\; R\rightarrow \infty .
\eeq
Letting $R\rightarrow \infty$ in \eqref{(10)}, and using \eqref{(11)}, \eqref{(12)}, we obtain
\beq\label{(13)}
\int_{\Bbb R} G(h(x))dx\geq 0
\eeq
provided
\beq\label{(14)}
\int_{C^L_R}G(h(z))ds\rightarrow  \int_{\Bbb R}G(h(x))dx \qquad \text{as}\; R\rightarrow \infty .
\eeq

To show \eqref{(14)}, using parametric equations, the integral $\int_{C^L_R}G(h(z))ds$ is equal to
\beq\label{(15)}
\int_{-R'\sqrt{2}/2}^{R'\sqrt{2}/2}G\bigg(h\bigg(x+iR-iR'\sqrt{1-\frac{x^2}{R'^2}}\bigg)\bigg)\frac{dx}{\sqrt{1-\frac{x^2}{R'^2}}}.
\eeq
In view of \eqref{(9)}, for all $R>100$, the integrand in \eqref{(15)} is bounded by the integrable function $C_f(1+|x|)^{-p}$ since $\sqrt{1-\frac{x^2}{R'^2}} $ is bounded from below by $\sqrt{1/2}$ in the range of integration. Then the Lebesgue dominated convergence theorem gives that \eqref{(15)} converges to
\beq\label{(40)}
\int_{\Bbb R}G(h(x))dx
\eeq
as $R\rightarrow \infty$.

Then replace $z$ with $h(x)=f(x)+iH(f)(x)$ in \eqref{(5)} and integrate \eqref{(5)} with respect to $x$, we get
\beq\label{(41)}
\int_{\Bbb R}|af(x)+bH(f)(x)|^pdx\leq B_p\int_{\Bbb R}|f(x)|^pdx-\int_{\Bbb R}G(h(x))dx .
\eeq
So by \eqref{(13)} we obtain
\beq\label{(42)}
\|(aI+bH)f\|^p_{L^p(\Bbb R)}\leq B_p\|f\|^p_{L^p(\Bbb R)}.
\eeq

\section{The sharpness of the constant $B_p$}\label{Sec5}

To deduce that the constant $B_p$ is sharp, we need to show
\beq\label{(19)}
\|aI+bH\|^p_{L^p({\Bbb R})\rightarrow L^p({\Bbb R})}\geq B_p.
\eeq
The proof of \eqref{(19)} relies on finding  suitable analytic functions in $H^p$ of the upper half space that will serve as
approximate extremals. Unlike the case of the circle,
where the   functions $\big((1+z)/(1-z)\big)^{1/p-\epsilon}$ in $H^p$ of the unit disc
serve  this purpose for all $1<p<\infty$ (see \cite{HKV}) as $\epsilon\downarrow 0$, we need to
consider the    cases $p<2$ and $p>2$ separately.\smallskip

\emph{Case 1:} $1<p<2$. Recall the analytic function used in \cite{GK} (also used in \cite{Pi}),
$$
F(z)=(z+1)^{-1}\bigg(i\frac{z+1}{z-1} \bigg)^{2\gamma /\pi}
$$
on the upper half plane. If $1<p<2$ and $\pi/2p'<\gamma<\pi/2p$, where $p'=p/(p-1)$, then $F(z)$ belongs to $H^p$ (the Hardy Spaces) in the upper half plane. Let
$$
f_{\gamma}(x)=\frac{1}{x+1}\bigg(\frac{|x+1|}{|x-1|} \bigg)^{2\gamma /\pi}\cos \gamma ,
$$
then we have
$$
F(x+i0)=f_{\gamma}(x)+i\begin{cases}
 \frac{1}{x+1}\big(\frac{|x+1|}{|x-1|} \big)^{2\gamma /\pi}\sin \gamma \quad\quad \text{when} \;|x|>1,\\
 -\frac{1}{x+1}\big(\frac{|x+1|}{|x-1|} \big)^{2\gamma /\pi}\sin \gamma \quad\; \text{when} \;|x|<1
 \end{cases}
$$
and since this is equal to the boundary values of a holomorphic function on the upper half plane, it follows that
$$
 H(f_{\gamma})(x)=\begin{cases}
 (\tan \gamma) f_{\gamma}(x) \quad\quad \text{when} \;|x|>1,\\
 -(\tan \gamma) f_{\gamma}(x) \quad\; \text{when} \;|x|<1,
 \end{cases}
$$
 So consider a function of the form $g_{\gamma}=\alpha f_{\gamma}+\beta H(f_{\gamma})$, where $\alpha ,\beta \in {\Bbb R}$. Notice that $H(g_{\gamma})=\alpha H(f_{\gamma})-\beta f_{\gamma}$, and the function $(|x-1|^{-\frac{2\gamma}{\pi}}|x+1|^{\frac{2\gamma}{\pi}-1})^p $ is integrable over the entire line since $\pi/2p'<\gamma<\pi/2p$, so for fixed $\alpha ,\beta$ we have
$$
 \begin{array}{cl}
 &\dfrac{\|(aI+bH)g_{\gamma}\|^p_{L^p({\Bbb R})}}{\|g_{\gamma}\|^p_{L^p({\Bbb R})}}
 =\dfrac{\int_{\Bbb R} |(a\alpha-b\beta)f_{\gamma}+(a\beta+b\alpha)H(f_{\gamma})|^pdx}{\int_{\Bbb R}|\alpha f_{\gamma}+\beta H(f_{\gamma})|^pdx} \\
 =&\dfrac{|(a\alpha-b\beta)+(a\beta+b\alpha)\tan \gamma|^pA_{\gamma}+|(a\alpha-b\beta)-(a\beta+b\alpha)\tan \gamma|^pB_{\gamma}}{|\alpha+\beta \tan \gamma|^pA_{\gamma}+|\alpha-\beta \tan \gamma|^pB_{\gamma}}
 \end{array}
$$
 where $A_{\gamma}=\int_{|x|>1}|f_{\gamma}(x)|^pdx, B_{\gamma}=\int_{|x|<1}|f_{\gamma}(x)|^pdx$. It is easy to get $A_{\gamma}\geq B_{\gamma}$, so
 \beq\label{(16)}
 \begin{array}{cl}
 &\dfrac{\|(aI+bH)g_{\gamma}\|^p_{L^p({\Bbb R})}}{\|g_{\gamma}\|^p_{L^p({\Bbb R})}}\\
 \geq &\dfrac{B_{\gamma}}{A_{\gamma}} \dfrac{|(a\alpha-b\beta)+(a\beta+b\alpha)\tan \gamma|^p+|(a\alpha-b\beta)-(a\beta+b\alpha)\tan \gamma|^p}{|\alpha+\beta \tan \gamma|^p+|\alpha-\beta \tan \gamma|^p},
 \end{array}
 \eeq
 and
 \beq\label{(17)}
 \begin{array}{cl}
 &\dfrac{\|(aI+bH)g_{\gamma}\|^p_{L^p({\Bbb R})}}{\|g_{\gamma}\|^p_{L^p({\Bbb R})}}\\
 \leq &\dfrac{A_{\gamma}}{B_{\gamma}} \dfrac{|(a\alpha-b\beta)+(a\beta+b\alpha)\tan \gamma|^p+|(a\alpha-b\beta)-(a\beta+b\alpha)\tan \gamma|^p}{|\alpha+\beta \tan \gamma|^p+|\alpha-\beta \tan \gamma|^p}.
 \end{array}
 \eeq
 Now we argue that
   \beq\label{(18)}
 \lim_{\gamma\rightarrow \frac{\pi}{2p}}\frac{A_{\gamma}}{B_{\gamma}}=1.
  \eeq
 In fact, by the second mean value theorem for definite integrals, there exists $\varepsilon \in (\delta ,1)$ where $0<\delta<1$ so that
$$
 \frac{\dint_{\delta}^1|x|^{p-2}\frac{|x+1|^{\frac{2\gamma p}{\pi}-p}}{|x-1|^{\frac{2\gamma p}{\pi}}}dx}{\dint_{\delta}^1\frac{|x+1|^{\frac{2\gamma p}{\pi}-p}}{|x-1|^{\frac{2\gamma p}{\pi}}}dx}
 =\frac{\frac{1}{\delta ^{2-p}}\dint_{\delta}^{\varepsilon}\frac{|x+1|^{\frac{2\gamma p}{\pi}-p}}{|x-1|^{\frac{2\gamma p}{\pi}}}dx+\dint_{\varepsilon}^1\frac{|x+1|^{\frac{2\gamma p}{\pi}-p}}{|x-1|^{\frac{2\gamma p}{\pi}}}dx}{\dint_{\delta}^{\varepsilon}\frac{|x+1|^{\frac{2\gamma p}{\pi}-p}}{|x-1|^{\frac{2\gamma p}{\pi}}}dx+\dint_{\varepsilon}^1\frac{|x+1|^{\frac{2\gamma p}{\pi}-p}}{|x-1|^{\frac{2\gamma p}{\pi}}}dx}.
$$
 Since $\int_{\varepsilon}^1|x+1|^{\frac{2\gamma p}{\pi}-p}|x-1|^{-\frac{2\gamma p}{\pi}}dx \rightarrow \infty$ as $\gamma\rightarrow \frac{\pi}{2p}$, we get
$$
 \lim_{\gamma\rightarrow \frac{\pi}{2p}} \frac{\dint_{\delta}^1|x|^{p-2}|x+1|^{\frac{2\gamma p}{\pi}-p}|x-1|^{-\frac{2\gamma p}{\pi}}dx}{\dint_{\delta}^1|x+1|^{\frac{2\gamma p}{\pi}-p}|x-1|^{-\frac{2\gamma p}{\pi}}dx}=1 .
$$
 Clearly this implies \eqref{(18)}. Combining   \eqref{(16)} with \eqref{(17)} we obtain
$$
 \begin{array}{cl}
 &\|aI+bH\|^p_{L^p({\Bbb R})\rightarrow L^p({\Bbb R})}\\
 \geq &\max\limits_{\alpha,\beta \in {\Bbb R}} \bigg(\dfrac{|(a\alpha-b\beta)+(a\beta+b\alpha)\tan \gamma '|^p+|(a\alpha-b\beta)-(a\beta+b\alpha)\tan \gamma '|^p}{|\alpha+\beta \tan \gamma '|^p+|\alpha-\beta \tan \gamma '|^p}\bigg)^{\frac{1}{p}},
 \end{array}
$$
 where $\gamma '=\frac{\pi}{2p}$. Letting $x=\alpha/\beta$ in \eqref{(3)}, we see that \eqref{(19)} holds, therefore the constant $B_p$ is sharp for $1<p<2$.

 \emph{Case 2:} $2<p<\infty$. In this case, the function $(|x-1|^{-\frac{2\gamma}{\pi}}|x+1|^{\frac{2\gamma}{\pi}-1})^p $ used in Case 1 fails to be integrable over the entire line. So we  consider   the following analytic function:
$$
 F(z)= \big(i(z^2-1)\big)^{-\frac{2\gamma}{\pi}},
 $$
 which belongs to $H^p$ in the upper half plane when $2<p<\infty$ and $\pi/4p<\gamma <\pi/2p$. Let
$$
f_{\gamma}(x)=|x+1|^{-\frac{2\gamma}{\pi}}|x-1|^{-\frac{2\gamma}{\pi}}\cos \gamma ,
$$
then we have
$$
F(x+i0)=f_{\gamma}(x)+i\begin{cases}
 -|x+1|^{-\frac{2\gamma}{\pi}}|x-1|^{-\frac{2\gamma}{\pi}} \sin \gamma \quad\quad \text{when} \;|x|>1,\\
 \quad|x+1|^{-\frac{2\gamma}{\pi}}|x-1|^{-\frac{2\gamma}{\pi}} \sin \gamma \quad\;\:\: \text{when} \;|x|<1.
 \end{cases}
$$
 It follows that
 $$
 H(f_{\gamma})(x)=\begin{cases}
 (\tan \gamma) f_{\gamma}(x) \quad\quad \text{when} \;|x|<1,\\
 -(\tan \gamma) f_{\gamma}(x) \quad\; \text{when} \;|x|>1,
 \end{cases}
$$
 Consider the function $g_{\gamma}=\alpha f_{\gamma}+\beta H(f_{\gamma})$, where $\alpha ,\beta \in {\Bbb R}$. Notice that the function $(|x-1|^{-\frac{2\gamma}{\pi}}|x+1|^{-\frac{2\gamma}{\pi}})^p $ is integrable over the entire line since $\pi/4p<\gamma <\pi/2p$, so for fixed $\alpha ,\beta$ we have
$$
 \begin{array}{cl}
 &\dfrac{\|(aI+bH)g_{\gamma}\|^p_{L^p({\Bbb R})}}{\|g_{\gamma}\|^p_{L^p({\Bbb R})}}\\
 =&\dfrac{|(a\alpha-b\beta)+(a\beta+b\alpha)\tan \gamma|^pA_{\gamma}+|(a\alpha-b\beta)-(a\beta+b\alpha)\tan \gamma|^pB_{\gamma}}{|\alpha+\beta \tan \gamma|^pA_{\gamma}+|\alpha-\beta \tan \gamma|^pB_{\gamma}}
 \end{array}
$$
 where $A_{\gamma}=\int_{|x|<1}|f_{\gamma}(x)|^pdx, B_{\gamma}=\int_{|x|>1}|f_{\gamma}(x)|^pdx$. It is easy to see $A_{\gamma}\leq B_{\gamma}$, so
 \beq\label{(n5)}
 \begin{array}{cl}
 &\dfrac{\|(aI+bH)g_{\gamma}\|^p_{L^p({\Bbb R})}}{\|g_{\gamma}\|^p_{L^p({\Bbb R})}}\\
 \leq &\dfrac{B_{\gamma}}{A_{\gamma}} \dfrac{|(a\alpha-b\beta)+(a\beta+b\alpha)\tan \gamma|^p+|(a\alpha-b\beta)-(a\beta+b\alpha)\tan \gamma|^p}{|\alpha+\beta \tan \gamma|^p+|\alpha-\beta \tan \gamma|^p},
 \end{array}
 \eeq
 and
 \beq\label{(n6)}
 \begin{array}{cl}
 &\dfrac{\|(aI+bH)g_{\gamma}\|^p_{L^p({\Bbb R})}}{\|g_{\gamma}\|^p_{L^p({\Bbb R})}}\\
 \geq &\dfrac{A_{\gamma}}{B_{\gamma}} \dfrac{|(a\alpha-b\beta)+(a\beta+b\alpha)\tan \gamma|^p+|(a\alpha-b\beta)-(a\beta+b\alpha)\tan \gamma|^p}{|\alpha+\beta \tan \gamma|^p+|\alpha-\beta \tan \gamma|^p}.
 \end{array}
 \eeq
 By the second mean value theorem for definite integrals, there exists $\varepsilon \in (\delta ,1)$ where $0<\delta<1$ so that
 \beq
 \begin{array}{cl}
 &\frac{\dint_{\delta}^1|x|^{\frac{4\gamma p}{\pi}-2}|x+1|^{-\frac{2\gamma p}{\pi}}|x-1|^{-\frac{2\gamma p}{\pi}}dx}{\dint_{\delta}^1|x+1|^{-\frac{2\gamma p}{\pi}}|x-1|^{-\frac{2\gamma p}{\pi}}dx}\notag \\
 =&\frac{\delta^{\frac{4\gamma p}{\pi}-2}\dint_{\delta}^{\varepsilon}|x+1|^{-\frac{2\gamma p}{\pi}}|x-1|^{-\frac{2\gamma p}{\pi}}dx+\dint_{\varepsilon}^1|x+1|^{-\frac{2\gamma p}{\pi}}|x-1|^{-\frac{2\gamma p}{\pi}}dx}{\dint_{\delta}^{\varepsilon}|x+1|^{-\frac{2\gamma p}{\pi}}|x-1|^{-\frac{2\gamma p}{\pi}}dx+\dint_{\varepsilon}^1|x+1|^{-\frac{2\gamma p}{\pi}}|x-1|^{-\frac{2\gamma p}{\pi}}dx}. \notag
 \end{array}
 \eeq
 Since $\int_{\varepsilon}^1|x+1|^{-\frac{2\gamma p}{\pi}}|x-1|^{-\frac{2\gamma p}{\pi}}dx \rightarrow \infty$ as $\gamma\rightarrow \frac{\pi}{2p}$, we have
$$
 \lim_{\gamma\rightarrow \frac{\pi}{2p}} \frac{\dint_{\delta}^1|x|^{\frac{4\gamma p}{\pi}-2}|x+1|^{-\frac{2\gamma p}{\pi}}|x-1|^{-\frac{2\gamma p}{\pi}}dx}{\dint_{\delta}^1|x+1|^{-\frac{2\gamma p}{\pi}}|x-1|^{-\frac{2\gamma p}{\pi}}dx}=1.
$$
 This implies
 $$
 \lim_{\gamma\rightarrow \frac{\pi}{2p}}\frac{B_{\gamma}}{A_{\gamma}}=1.
$$
 Combining   \eqref{(n5)} and \eqref{(n6)} we obtain
 $$
 \begin{array}{cl}
 &\|aI+bH\|^p_{L^p({\Bbb R})\rightarrow L^p({\Bbb R})}\\
 \geq &\max\limits_{\alpha,\beta \in {\Bbb R}} \bigg(\dfrac{|(a\alpha-b\beta)+(a\beta+b\alpha)\tan \frac{\pi}{2p}|^p+|(a\alpha-b\beta)-(a\beta+b\alpha)\tan \frac{\pi}{2p}|^p}{|\alpha+\beta \tan \frac{\pi}{2p}|^p+|\alpha-\beta \tan \frac{\pi}{2p}|^p}\bigg)^{\frac{1}{p}}.
 \end{array}
$$
 Letting $x=\alpha/\beta$ in \eqref{(3)}, so \eqref{(19)} holds, therefore the constant $B_p$ is sharp for $2<p<\infty$.\\ \qed


\begin{thebibliography}{120}

\bibitem{decarli-laeng} L. De Carli, E. Laeng,
\emph{Sharp $L^p$ estimates for the segment multiplier}, Collect. Math. {\bf 51} (2000),   309--326.

\bibitem{Gamelin}
T. W. Gamelin, \emph{Uniform Algebras and Jensen Measures},  London Math. Soc. Lecture Note Series, Vol. 32, Cambridge Univ. Press, Cambridge New York, 1978.


\bibitem{G} L. Grafakos, \emph{Classical Fourier Analysis}, 3r ed.,   GTM 249, Springer
New York,  2014.

\bibitem{G2} L. Grafakos, \emph{Best bounds for the Hilbert transform on $L^p(\mathbb R^1)$}, Math. Res. Let.   {\bf 4} (1997), 469--471.

\bibitem{GK} T. Gokhberg, N. Y. Krupnik, \emph{Norm of the Hilbert transformation in the $L^p$ space},
Funktsional. Analiz i Ego Prilozhen  {\bf 2} (1968), 91--92.

\bibitem{GS} L. Grafakos, T. Savage, \emph{Best bounds for the Hilbert transform on $L^p(\mathbb R^1)$; A corrigendum}, Math. Ress Let,  {\bf 22} (2015),   1333--1335.

\bibitem{HKV} B. Hollenbeck, N. J. Kalton, I. E. Verbitsky, \emph{Best constants for some operators associated with the Fourier and Hilbert transforms}, Studia Math.  {\bf 157} (2003), 237--278.

\bibitem{HV} B. Hollenbeck,  I. E. Verbitsky,   \emph{Best constants for the Riesz projection}, J.  Funct. Anal. {\bf 175} (2000),  370--392.

\bibitem{Laeng} E. Laeng, \emph{Remarks on the Hilbert transform and on some families of multiplier operators related to it}, Collect. Math.  {\bf 58} (2007),   25--44.


\bibitem{L} B. Ya. Levin, \emph{Lectures on entire functions}, Transl. Math. Monogr. 150, Amer. Math. Soc., Providence, RI, 1996.

\bibitem{MZ} J. Marcinkiewicz, A. Zygmund, \emph{Quelques in\'egalit\'es pour les op\'erations lin\'eaires},  Fund. Math. {\bf 32} (1939), 112--121.

\bibitem{Pi} S. K. Pichorides, \emph{On the best values of the constants in the theorems of M. Riesz, Zygmund and Kolmogorov}, Studia Math., {\bf 44} (1972), 165--179.

\bibitem{Z} A. Zygmund, \emph{Trigonometric Series}, Vol 2, Cambridge Univ. Press, London, UK, 1968.




\end {thebibliography}
\end{document}